\documentstyle[titlepage,twoside,12pt]{article}
\begin{document}

\pagenumbering{arabic}
\setcounter{page}{1}

\pagenumbering{arabic}

{\large \bf How to Solve Smooth Nonlinear PDEs in Algebras of Generalized Functions with Dense
Singularities} \\ \\

{\it Elem\'{e}r E Rosinger \\ Department of Mathematics \\ University of Pretoria \\ Pretoria, 0002 South Africa \\
e-mail : eerosinger@hotmail.com} \\ \\

{\bf Abstract} \\

As a significant strengthening of properties of earlier algebras of generalized functions,
here are presented classes of such algebras which can deal with {\it dense} singularities. In
fact, the cardinal of the set of singular points can be {\it larger} than that of the
nonsingular points. This large class of singularities allows the solution of large classes of
smooth nonlinear PDEs. This in certain ways overcomes the celebrated 1957 H. Lewy
impossibility result. \\

Note : This paper is an updated version of the paper with the same title published in
Applicable Analysis, vol. 78, 2001, pp. 355-378. \\ \\

{\bf 1. The Algebraic-Differential and the Order Completion Solution Methods} \\

As an illustration of the interest and power of the recently introduced differential algebras
of generalized functions with {\it dense} singularities, see Rosinger [14-17] and Mallios \&
Rosinger [2], we shall present here a general, that is, {\it type independent} existence
result which gives solutions for rather arbitrary smooth nonlinear PDEs in these algebras.
This existence of solutions result has an interest of its own since, among others, it helps
overcoming in part the celebrated 1957 impossibility of H. Lewy, related to the solution of
smooth linear PDEs, and it does so in the general smooth {\it nonlinear} case, see Comments in
section 8. \\
Here it can be mentioned that the Lewy impossibility had for the first time been completely
overcome in Oberguggenberger \& Rosinger, and it was done so with a large margin, namely,
within very general conditions of {\it nonlinearity} and {\it lack of smoothness} on the PDEs
involved. Indeed, in the work cited, with the use of the order completion method introduced
there and applied to spaces of smooth functions, all continuous nonlinear PDEs were given
solutions which can be assimilated with usual {\it measurable} functions. More recently,
however, this universal regularity result has been improved by showing that all continuous
nonlinear PDEs have solutions which can be assimilated with Hausdorff continuous functions,
see Anguelov \& Rosinger. \\
Nevertheless, when compared with such existence results obtained by order completion, a
possible advantage of the existence result presented in Theorem 1 in the sequel is that,
although less general since it requires smoothness of the PDEs, it nevertheless has a certain
increased {\it universality} property with respect to the algebras in which the solutions
prove to exist. Indeed, these algebras appear to depend less on the PDEs solved, than is the
case with the spaces of functions delivered by the order completion method. \\
Another feature of possible interest of the existence result presented here is that it came
about recently by a certain interaction between two rather different nonlinear theories of
generalized solutions, namely, the older algebraic one, and the more recent one, based on
order completion in Oberguggenberger \& Rosinger. Relevant references can be found in the
Comments in section 8. \\
And now, to the details of the type independent and rather universal kind of existence result
for smooth nonlinear PDEs. \\ \\

{\bf 2. The General Class of Smooth Nonlinear PDEs} \\

The smooth nonlinear PDEs considered are of the general form

\bigskip
(2.1) \quad $ T(x,D) U(x) ~=~ f(x),~~~ x \in X $

\medskip
where $X$ is a nonvoid, possibly unbounded open subset of ${\bf R}^n$, $f$ is a smooth
function on $X$, $U$ is the unknown function, while the smooth nonlinear partial differential
operator in the left hand term is given by

\bigskip
(2.2) \quad $ T(x,D) U(x) ~=~ F (x, U(x), ~.~.~.~ , D^p U(x), ~.~.~.~ ),~~~ x \in X $

\medskip
where $F$ is any function jointly smooth in all its arguments, while $p \in {\bf N}^n$,
$| p | \leq m$, for an arbitrary given $m \in {\bf N}$. All these smooth functions, as well as
the generalized functions which appear in the sequel are considered for simplicity real valued.
However, the extension to the case of finite dimensional vector valued functions, and in
particular, complex valued ones, is rather immediate. \\

{\bf Range Conditions}. Let us consider for $x \in X$ the subset of real numbers

\bigskip
(2.3) \quad $ R_x ~=~ \{~ F(x, \xi_0, ~.~.~.~ , \xi_p, ~.~.~.~ ) ~|~ \xi_p \in {\bf R},~ p
                                \in {\bf N}^n,~ |p| \leq m ~\} $

\medskip
which is the {\it range} in ${\bf R}$ of $F(x, ~.~.~.~ )$, for every fixed $x \in X$. Since
$F$ is jointly continuous in its arguments, the set $R_x$ must be an interval in ${\bf R}$.
And in the case of nontrivial PDEs,  if (2.1) is linear, as well as in most of the nonlinear
cases of applicative interest, it will turn out that we have

\bigskip
(2.4) \quad $ R_x ~=~ {\bf R},~~~ x \in X $

\medskip
Now given any $x_0 \in X$,  an obvious {\it necessary} condition for the existence of a {\it
classical} solution of (2.1) in a neighbourhood of $x_0$ is that

\bigskip
(2.5) \quad $ f(x_0) \in R_{x_0} $

\medskip
In the sequel, related to the smooth nonlinear PDEs (2.1) we shall also consider certain
variants of the {\it range condition} (2.5). One of them will be the somewhat stronger
condition

\bigskip
(2.6) \quad $ f(x) \in ~\mbox{int}~ R_x,~~~ x \in X $

\medskip
Other variants will be derived from the significantly weaker version of (2.5), namely

\bigskip
(2.7) \quad $ f(x) \in R_x,~~~ x \in X_{sol} $

\medskip
where $X_{sol} \subseteq X$ is a suitably given nonvoid subset of $X$. \\
Clearly, whenever (2.4) holds, then condition (2.6), and in particular (2.7) will be satisfied.
It follows that in the nontrivial case, for linear smooth PDEs, as well as for most of the
nonlinear smooth PDEs of applicative interest, both conditions (2.6) and (2.7) will
automatically hold, as a consequence of (2.4). \\

In this way, as will be later argued, the existence result in Theorem 1 in section 4 does
indeed help to a good extent to overcome the H. Lewy impossibility. \\ \\

{\bf 3. The Differential Algebras of Generalized Functions with Dense Singularities} \\

The differential algebras of generalized functions in which we shall find solutions $U$ for
(2.1) are of the form

\bigskip
(3.1) \quad $ B (X) ~=~ ( {\cal C}^\infty (X) )^{\bf N} / {\cal J} (X) $

\medskip
with suitably constructed ideals ${\cal J} (X)$ in $( {\cal C}^\infty (X) )^{\bf N}$. \\
The construction of similar ideals, introduced and employed recently in Rosinger [14-16] and
Mallios \& Rosinger [2], has actually had a different motivation. Namely, one of them was
coming from the need to create a suitable mathematical framework for the so called 'space-time
foam' singularity structures which are of interest in General Relativity, structures supposed
to allow {\it dense} singularities in the four dimensional Einsteinian manifolds. The other
motivation is related to a similar need for highly singular structures in Quantum Gravity.
For further details, see Comments in section 8. \\

The basic idea in this respect is to construct ideals which, in the sense specified in the
sequel, can handle through algebraic-differential means the {\it largest} possible {\it
singularity} subsets $\Sigma$ in $X$. \\
Until recently, the largest such singularity subsets $\Sigma$ were supposed to be closed and
nowhere dense in $X$, see Rosinger [3-13] and Mallios \& Rosinger [1]. That itself was not
trivial, since such subsets can have arbitrary large positive Lebesgue measure, Oxtoby.
Furthermore, with the help of such ideals and the corresponding differential algebras of
generalized functions, one could already obtain type independent nonlinear existence results
such as the {\it global} solution of all analytic nonlinear PDEs, see for instance Rosinger
[7-9]. \\

{\bf Dense Singularities}. This time however, we shall be able to include in our
algebraic-differential approach {\it far larger} singularity subsets $\Sigma$ in $X$, namely,
all those which satisfy the rather mild condition that their complementary, that is, the set
of nonsingular points $X \setminus \Sigma$, is dense in X

\bigskip
(3.2) \quad $ X \setminus \Sigma ~~\mbox{is dense in}~ X $

\medskip
In particular, and as an {\it extreme} situation, it will be sufficient if $X \setminus
\Sigma$ is merely countable and dense in $X$, in which case the singularity set $\Sigma$
itself is dense in $X$, and on top of it, it is uncountable, thus, it has a larger cardinal
than the nonsingularity set $X \setminus \Sigma$, or to summarize

\bigskip
(3.3) \quad $ \begin{array}{l}
                             X \setminus \Sigma ~~\mbox{is countable and dense in}~ X \\
                             \\
                             \Sigma ~~\mbox{is uncountable and dense in}~ X \\
                             \\
                             \mbox{car}~ \Sigma ~>~ \mbox{car}~ X \setminus \Sigma
                       \end{array} $ \\

\bigskip
{\bf Families of Singularities}. What is further important is that we shall be able to deal
not only individually with such singularity sets $\Sigma$, but also with whole families
${\cal S}$ made up of them, provided that such families satisfy the two conditions

\bigskip
(3.4) \quad $ \begin{array}{l}
                         \forall~ \Sigma \in {\cal S} ~: \\
                         \\
                         ~~~ X \setminus \Sigma ~~\mbox{is dense in}~ X
                      \end{array} $

\medskip
which is but a repetition of (3.2), and

\bigskip
(3.5) \quad $ \begin{array}{l}
                         \forall~ \Sigma,~ \Sigma^\prime \in {\cal S} ~: \\
                         \\
                         \exists~ \Sigma^{\prime \prime} \in {\cal S} ~: \\
                         \\
                         ~~~ \Sigma \cup \Sigma^\prime \subseteq \Sigma^{\prime \prime}
                      \end{array} $ \\

\bigskip
{\bf Maximal Families of Singularities}. A direct application of Zorn's lemma, thus, of the
Axiom of Choice, gives for each singularity subset $\Sigma$ at least one {\it maximal} family
of singularity subsets ${\cal S}$ to which it belongs. In other words

\bigskip
(3.6) \quad $ \begin{array}{l}
                          \forall~ \Sigma \subset X ~\mbox{which satisfies}~ (3.2) ~: \\
                          \\
                          \exists~ \mbox{maximal}~ {\cal S} ~\mbox{satisfying}~ (3.4), (3.5) ~: \\
                          \\
                          ~~~ \Sigma \in {\cal S}
                       \end{array} $ \\

\bigskip
{\bf The Ideals}. And now, let us define the respective ideals. Given a singularity subset
$\Sigma \subset X$ which satisfies (3.2), we denote by

\bigskip
(3.7) \quad $ {\cal J}_\Sigma (X) $

\medskip
the {\it ideal} in $( {\cal C}^\infty (X) )^{\bf N}$ of all the sequences of smooth functions
$w = ( w_\nu ~|~ \nu \in {\bf N} )$ which {\it outside} of the singularity set $\Sigma$
satisfy the {\it asymptotic vanishing} condition

\bigskip
(3.8) \quad $ \begin{array}{l}
                          \forall~ x \in X \setminus \Sigma,~ l \in {\bf N} ~: \\
                          \\
                          \exists~ \nu \in {\bf N} ~: \\
                          \\
                          \forall~ \mu \in {\bf N},~ \mu \geq \nu,~ p \in {\bf N}^n,~ |p| \leq l ~: \\
                          \\
                          ~~~D^p w_\mu (x) ~=~ 0
                       \end{array} $

\medskip
Here we should note that the above asymptotic vanishing condition is of {\it finite} type in
the sense that at points $x$ outside of the singularity set $\Sigma$ only a finite - even if
arbitrary large - number of partial derivatives are required to vanish. However, by
considering different vanishing conditions for the same singularites one can, among others,
note the versatility of the algebraic approach in their study, see details in Rosinger
[14-17] and Mallios \& Rosinger [2]. \\ \\

{\bf Example 1.1} \\

Let $x_\nu \in X$, with $\nu \in {\bf N}$, be any given dense sequence in $X$, and let us
denote by $\Sigma$ its complementary, that is

$$ \Sigma ~=~ X \setminus \{~ x_\nu ~|~ \nu \in {\bf N} ~\} $$

\medskip
Then $\Sigma$ clearly satisfies (3.3). \\

Further, let any $l _\nu \in {\bf N}$, with $\nu \in {\bf N}$, be such that
$\lim_{\nu \to \infty}~ l_\nu = \infty$. \\

Now we define the sequence of polynomials $w = ( w_\nu ~|~ \nu \in {\bf N} )$ in the variable
$x \in X$, according to

$$ w_\nu (x) = ( x - x_0 )^{l_\nu} ~.~.~.~ ( x - x_\nu )^{l_\nu},~~ \nu \in {\bf N},~
                       x \in X $$

\medskip
Then in view of (3.8), it is easy to see that

$$ w \in {\cal J}^f_\Sigma (X) $$

\hfill $\diamondsuit \diamondsuit \diamondsuit$ \\

Based on the above, if we are given a family ${\cal S}$ of singularity subsets
$\Sigma \subset X$, family which satisfies (3.4), (3.5), then we associate with it the {\it
ideal} in $( {\cal C}^\infty (X) )^{\bf N}$ defined by

\bigskip
(3.9) \quad $ {\cal J}_{\cal S} (X) ~=~ \bigcup_{\Sigma \in {\cal S}}~ {\cal J}_\Sigma (X) $

\medskip
Corresponding to the ideals (3.7), (3.9), we define the algebras of generalized functions, see
(3.1)

\bigskip
(3.10) \quad $ \begin{array}{l}
                B_\Sigma (X) ~=~  ( {\cal C}^\infty (X) )^{\bf N} / {\cal J}_\Sigma (X) \\ \\
                B_{\cal S} (X) ~=~  ( {\cal C}^\infty (X) )^{\bf N} / {\cal J}_{\cal S} (X)
               \end{array} $

\medskip
Obviously if (3.2) holds for a singularity subset $\Sigma \subset X$ then the family
${\cal S} = \{~ \Sigma~\}$ consisting of this single singularity subset will satisfy (3.4),
(3.5). In this way, the algebras $B_\Sigma (X)$ are particular cases of the algebras
$B_{\cal S} (X)$. \\

{\bf Properties of the Algebras}. Some of the most important properties of these algebras
result from the following two properties of the respective ideals

\bigskip
(3.11) \quad $ {\cal J}_\Sigma (X) \bigcap {\cal U} (X) ~=~ \{~ 0 ~\},~~~
                         {\cal J}_{\cal S} (X) \bigcap {\cal U} (X) ~=~ \{~ 0 ~\} $

\medskip
and

\bigskip
(3.12) \quad $ D^p {\cal J}_\Sigma (X) \subseteq {\cal J}_\Sigma (X),~~~
                D^p {\cal J}_{\cal S} (X) \subseteq {\cal J}_{\cal S} (X),~~~ p \in
                           {\bf N}^n $

\medskip
Here in (3.11), ${\cal U} (X)$ is the {\it diagonal} in $( {\cal C}^\infty (X) )^{\bf N}$,
that is, it is the subalgebra of all the sequences $u(\psi) = (\psi, \psi, \psi, ~.~.~.~ )$,
with $\psi \in {\cal C}^\infty (X)$. In this way (3.11) means that the respective ideals are
{\it off diagonal}. As we shall see later, this is one of the most important properties of the
ideals used in the construction of algebras of generalized functions. \\
As far as (3.12) is concerned, this expresses the fact that the respective ideals are
invariant under arbitrary partial derivation. \\
Now it follows that (3.11) is the {\it necessary and sufficient} condition for the existence
of the {\it algebra embeddings}

\bigskip
(3.13) \quad $ \begin{array}{l}
                 {\cal C}^\infty (X) \ni \psi \longmapsto u(\psi) + {\cal J}_\Sigma (X)
                                   \in B_\Sigma (X) \\ \\
                 {\cal C}^\infty (X) \ni \psi \longmapsto u(\psi) + {\cal J}_{\cal S} (X)
                               \in B_{\cal S} (X)
               \end{array} $

\medskip
On the other hand, (3.12) turns the respective algebras of generalized functions into {\it
differential algebras}, according to the partial derivatives of their generalized functions
which can be defined as follows

\bigskip
(3.14) \quad $ \begin{array}{l}
               B_\Sigma (X) \ni U = s + {\cal J}_\Sigma (X) \longmapsto D^p U =
                         D^p s + {\cal J}_\Sigma (X) \in B_\Sigma (X) \\
                            \\
               B_{\cal S} (X) \ni U = s + {\cal J}_{\cal S} (X) \longmapsto D^p U =
                         D^p s + {\cal J}_{\cal S} (X) \in B_{\cal S} (X)
              \end{array} $

\medskip
where for $s = ( s_\nu ~|~ \nu \in {\bf N} ) \in ( {\cal C}^\infty (X) )^{\bf N}$ and
$p \in {\bf N}^n$, we denoted $D^p s = ( D^p s_\nu ~|~ \nu \in {\bf N} )$, that is, the
termwise partial derivation of the sequence of smooth functions $s$. \\
In this way, and in view of (3.2), we obtain the {\it differential algebras of generalized
functions with dense singularities} $B_\Sigma (X)$ and $B_{\cal S} (X)$. Moreover, the algebra
embeddings (3.13) are in fact {\it embeddings of differential algebras}. \\

Finally, as seen in Rosinger [1-17], the off diagonality properties (3.11) also imply that
these algebras contain the Schwartz distributions as vector subspaces, namely

\bigskip
(3.15) \quad $ {\cal D}^\prime (X) \subset B_\Sigma (X),~~~ {\cal D}^\prime (X) \subset
                              B_{\cal S} (X) $ \\

\bigskip
{\bf Increased Regularity}. Let us take two singularity subsets $\Sigma,~ \Sigma^\prime
\subset X$ which satisfy (3.2), and two families of singularities ${\cal S},~ {\cal S}^\prime$
for which (3.4) and (3.5) hold. It is easy to see that in such a case we shall have the
implications

\bigskip
(3.16) \quad $ \begin{array}{l}
                           \Sigma \subseteq \Sigma^\prime ~\Longrightarrow~
                            {\cal J}_\Sigma (X) \subseteq {\cal J}_{\Sigma^\prime} (X) \\ \\
                           {\cal S} \subseteq {\cal S}^\prime ~\Longrightarrow~
                          {\cal J}_{\cal S} (X) \subseteq {\cal J}_{{\cal S}^\prime} (X) \\ \\
                           \Sigma \in {\cal S} ~\Longrightarrow~
                             {\cal J}_\Sigma (X) \subseteq {\cal J}_{\cal S} (X)
                \end{array} $

\medskip
As a consequence, when the premises in these implications hold, we obtain the {\it canonical
surjective algebra homomorphisms}

\bigskip
(3.17) \quad $ \begin{array}{l}
                             B_\Sigma (X) \ni s + {\cal J}_\Sigma (X) \longmapsto
                                   s + {\cal J}_{\Sigma^\prime} (X) \in
                                           B_{\Sigma^\prime} (X) \\
                             \\
                             B_{\cal S} (X) \ni s + {\cal J}_{\cal S} (X) \longmapsto
                            s + {\cal J}_{{\cal S}^\prime} (X) \in B_{{\cal S}^\prime} (X) \\
                             \\
                              B_\Sigma (X) \ni s + {\cal J}_\Sigma (X) \longmapsto
                                      s + {\cal J}_{\cal S} (X) \in B_{\cal S} (X)
                         \end{array} $

\medskip
and as seen in Rosinger [14-17], these canonical surjective algebra homomorphisms can in a
suitable sense be interpreted as giving algebras in their right hand terms whose elements are
{\it more regular} than those of the algebras in the left hand terms. \\

{\bf Partial Differential Operators in the Algebras}. Given any family ${\cal S}$ of
singularities which satisfies (3.4), (3.5), we can define the action, see for details
Rosinger [14-17]

\bigskip
(3.18) \quad $ T(x,D) : B_{\cal S} (X) ~\longrightarrow~ B_{\cal S} (X) $

\medskip
of the smooth nonlinear partial differential operator in (2.2) on the differential algebra of
generalized functions corresponding to ${\cal S}$ as follows. This action will be a natural
extension of the classical action

$$T(x,D) : {\cal C}^\infty (X) ~\longrightarrow~ {\cal C}^\infty (X) $$

\medskip
of the smooth nonlinear partial differential operator (2.2). First, and independently of $F$,
we simply collect together all the needed partial derivatives, that is, we define the mapping

$$ B_{\cal S} (X) \ni U ~\longrightarrow~ ( U, ~.~.~.~ , D^p U, ~.~.~.~ ) \in
                 ( B_{\cal S} (X) )^{m^*} $$

\medskip
where $p \in {\bf N}^n,~ |p| \leq m$, while $m^*$ is the number of real variable arguments in
$F$ minus $n$. Now, we compose the above mapping with the purely nonlinear and nondifferential
one defined by $F$ according to

$$ F : ( B_{\cal S} (X) )^{m^*} ~\longrightarrow~ B_{\cal S} (X) $$

\medskip
and in the following way. Given $U_p = s_p + {\cal J}_{\cal S} (X) \in B_{\cal S} (X),~ p \in
{\bf N}^n,~ |p| \leq m$, where $s_p = ( s_{p \nu} ~|~ \nu \in {\bf N} ) \in
( {\cal C}^\infty (X) )^{\bf N}$, we define

$$ F(U_0, ~.~.~.~ , U_p, ~.~.~.~ ) = t + {\cal J}_{\cal S} (X) \in B_{\cal S} (X) $$

\medskip
where $t = ( t_\nu ~|~ \nu \in {\bf N} ) \in  ( {\cal C}^\infty (X) )^{\bf N}$ is given by

$$ t_\nu (x) = F(x, s_{0 \nu} (x), ~.~.~.~ , s_{p \nu} (x), ~.~.~.~ ),~~ \nu \in {\bf N},~
                      x \in X $$

\medskip
Now it is easy to see that $T(x,D)$ in (3.18) commutes with the algebra embeddings (3.13) and
the classical action of (2.2), therefore, it is indeed a natural extension of the latter. \\
In this way, and in view of (3.13), the smooth nonlinear PDEs in (2.1) is well defined in each
differential algebra of generalized functions $B_{\cal S} (X)$, and for convenience, it can be
written in the form

\bigskip
(3.19) \quad $ T(x,D) U ~=~ f ~~~\mbox{in}~~ B_{\cal S} (X) $

\medskip
And since the algebras $B_\Sigma (X)$ are particular cases of $B_{\cal S} (X)$, the above
relations are thus well defined in the former algebras as well. \\ \\

{\bf 4. Solving the Smooth Nonlinear PDEs} \\

The differential algebras of generalized functions (3.10) which will contain the solutions $U$
of arbitrary smooth nonlinear PDEs (2.1) will be obtained by choosing various suitable ideals
(3.7). In this respect we have \\ \\

{\bf Theorem 1} \\

Given on a possibly unbounded nonvoid open subset $X \subseteq {\bf R}^n$ any smooth nonlinear
PDE in (2.1) which satisfies on a dense subset $X_{sol}$ of $X$ the stronger version of range
condition (2.7) presented in (5.20) below. \\
Then one can construct singularity subsets $\Sigma \subset X$, with $X \setminus \Sigma$ dense
in $X$, and $X \setminus \Sigma \subseteq X_{sol}$, together with corresponding generalized
functions

\bigskip
(4.1) \quad $ U = s + {\cal J}_\Sigma (X) \in B_\Sigma (X) $

\medskip
which in the sense of (3.19) are solutions in the differential algebra of generalized
functions  $B_\Sigma (X)$ of the PDE in (2.1), namely

\bigskip
(4.2) \quad $ T(x,D) U ~=~ f ~~~\mbox{in}~~ B_\Sigma (X) $

\hfill $\diamondsuit \diamondsuit \diamondsuit$ \\

The proof of Theorem 1 will follow from bringing together Corollary 2 in section 5, with the
above constructions related to the differential algebras of generalized functions with dense
singularities. \\
Next, in section 5, we shall shortly recall certain basic and rather unusual results from the
order completion method, results which recently suggested the vanishing property in
Corollary 2. Then in section 6, we shall construct the required solution algebras with dense
singularities needed in Theorem 1. \\

{\bf Categorial Properties}. It is easy to see that under the premises in the implications (3.16), the extension $T(x,D)$
in (3.18) of the classical action of the smooth nonlinear partial differential operator (2.2) commutes with the
canonical surjective algebra homomorpisms in (3.17). Therefore, under the respective premises, the solution $U$ in
(4.1) yields through the corresponding canonical mappings solutions in the algebras  $B_{\Sigma^\prime} (X)$,
$B_{{\cal S}^\prime} (X)$ and $B_{\cal S} (X)$. \\

A possible interest in these latter solutions is connected with the fact that, as mentioned, they can be interpreted as
being {\it more regular} than the original solution $U$, see also Rosinger [14-17]. \\ \\

{\bf 5. Some Basic Results from the Order Completion Method} \\

{\bf The Earlier Results}. We start by presenting step by step the basic approximation results
upon which the whole order completion solution method in Oberguggenberger \& Rosinger rests.
These results are rather unusual even if quite simple, and as such, they recently suggested a
certain extension which leads to the {\it vanishing} property of {\it error} sequences in
Corollary 2, upon which the existence result in Theorem 1 is based. \\
Next in Lemma 1, Proposition 1 and Corollary 1 we shall - as in Obergugenberger \& Rosinger -
assume the general case when $F$ and $f$ in (2.1), (2.2) are merely  jointly {\it continuous}
in their argument, and not necessarily smooth. On the other hand, we shall assume the stronger
range condition (2.6). The proofs of Lemma 1 and Proposition 1 can be found in
Oberguggenberger \& Rosinger. \\ \\

{\bf Lemma 1} \\

Given any nonvoid possibly unbounded open subset $X \subseteq {\bf R}^n$, then

\bigskip
(5.1) \quad $ \begin{array}{l}
                  \forall~ x_0 \in X,~ \epsilon > 0 ~: \\
                           \\
                  \exists~ \delta > 0,~ P ~\mbox{polynomial in}~ x \in {\bf R}^n ~: \\
                           \\
                  \forall~ x \in X ~: \\
                           \\
                   ~~~ \|x - x_0 \| \leq \delta ~~\Longrightarrow~~ f(x) - \epsilon ~\leq~
                                  T(x,D) P(x) ~\leq~ f(x)
                  \end{array} $

{~} \\ \\

{\bf Proposition 1} \\

Given any nonvoid possibly unbounded open subset $X \subseteq {\bf R}^n$, then

\bigskip
(5.2) \quad $ \begin{array}{l}
                \forall~ \epsilon > 0 ~: \\
                             \\
                \exists~ \Gamma_\epsilon \subset X,~ \Gamma_\epsilon ~
                         \mbox{closed and nowhere dense in}~ X, \\
                ~~~~~~~~~~~~~~~~~~~~~~~~~~~~~~~~~~~~~~~ U_\epsilon \in
                         {\cal C}^\infty ( X \setminus \Gamma_\epsilon ) ~: \\
                             \\
                       ~~~ f -\epsilon ~\leq~ T(x,D) U_\epsilon ~\leq~
                            f ~~\mbox{on}~ X \setminus \Gamma_\epsilon
               \end{array} $

{~} \\ \\
{\bf Remark 1} \\

1. We note that in (5.2) we can construct $\Gamma_\epsilon$ in such a way as to have the
additional property

\bigskip
(5.3) \quad $ \mbox{mes}~ \Gamma_\epsilon ~=~ 0 $

\medskip
2. By a similar argument, we can obtain a version of (5.2) in which we have the inequalities

\bigskip
(5.4) \quad $ f ~\leq~ T(x,D) U_\epsilon ~\leq~ f + \epsilon
                       ~~\mbox{on}~ X \setminus \Gamma_\epsilon $

\medskip
3. It is important to note that the presence of the {\it closed nowhere dense} subset
$\Gamma_\epsilon \subset X$ in (5.2) is in fact a quite natural, minimal and also unavoidable
type of lack of global regularity. Indeed, even in the particular case when both $T(x,D)$ and
$f$ are analytic, we still cannot expect to have classical solutions on the whole of the
domain of definition $X$ of the PDEs in (2.1), see Rosinger [14-17] for further details. \\ \\

{\bf Corollary 1} \\

Under the conditions in Proposition 1, we have the solution property

\bigskip
(5.5) \quad $ \begin{array}{l}
                    \forall~ x \in X ~: \\
                    \\
                    \exists~ \delta > 0 ~: \\
                    \\
                    \forall~ A \subset X \cap B(x, \delta),~ A ~\mbox{finite} ~: \\
                    \\
                    \exists~ U \in {\cal C}^\infty (X) ~: \\
                    \\
                    ~~~ T(x,D) U (y) = f(y),~~ y \in A
              \end{array} $

\medskip
where $B(x, \delta)$ denotes the open ball of radius $\delta$ in ${\bf R}^n$ around $x$.

\bigskip
{\bf Proof} \\
Given $x \in X$ and $\epsilon > 0$, we obtain from (5.1) a function $U_{-} \in
{\cal C}^\infty (X)$ such that

\bigskip
(5.6) \quad $ T(y,D) U_{-} (y) ~\leq~ f(y),~~ y \in X \cap B(x, \delta) $

\medskip
By a similar argument, we obtain a function $U_{+} \in {\cal C}^\infty (X)$ for which

\bigskip
(5.7) \quad $ f(y) ~\leq~ T(y,D) U_{+} (y),~~ y \in X \cap B(x, \delta) $

\medskip
Now for $\lambda \in [0,1]$ let us define the convex combination

\bigskip
(5.8) \quad $ U_\lambda ~=~ ( 1 - \lambda ) U_{-} + \lambda U_{+} \in {\cal C}^\infty (X) $

\medskip
and the continuous function $h : [0,1] \times X \longrightarrow {\bf R}$ by
$h(\lambda,y) ~=~ T(y,D) U_\lambda (y)$. Then
clearly $ h(0,y) ~\leq~ 0 ~\leq~ h(1,y),~~ y \in X \cap B(x, \delta) $, hence the continuity of $h$ results in

\bigskip
(5.9) \quad $ \begin{array}{l}
                    \forall~ y \in X \cap B(x, \delta) ~: \\
                    \\
                    \exists~ \lambda \in [0,1] ~: \\
                    \\
                    ~~~ T(y,D) U_\lambda (y) = f(y)
                \end{array} $

\medskip
since the above equality simply means that $h(\lambda,y) = 0$. \\
Let us now take any finite subset $A \subset X \cap B(x, \delta)$ and apply (5.8) to each $a \in A$, thus obtaining
$U_{\lambda_a} \in {\cal C}^\infty (X)$ which satisfies (2.1) at $a$. But as $A$ is finite, we can consider on $X$ a
partition of unity given by $\psi_a \in {\cal C}^\infty (X)$, with $a \in A$, such that

\bigskip
(5.10) \quad $ \begin{array}{l}
                   \psi_a = 1 ~~\mbox{on a neighbourhood of}~ a \\
                   \\
                   0 ~\leq~ \psi_a ~\leq~ 1 ~~\mbox{on}~ X \\
                   \\
                   \sum_{a \in A}~ \psi_a = 1 ~~\mbox{on}~ X
                \end{array} $

\medskip
Finally we define the function $U \in {\cal C}^\infty (X)$ by

\bigskip
(5.11) \quad $ U ~=~ \sum_{a \in A}~ \psi_a U_{\lambda_a} $

\medskip
and then the relations (5.8) - (5.11) will give (5.5). \\ \\

{\bf Remark 2} \\

One can note that the function $U$ in (5.11) which gives the result in Corollary 1 is simply
the convex combination, this time with variable coefficients, of the two functions $U_{-}$ and
$U_{+}$ in (5.6) and (5.7)  respectively, namely

$$ U ~=~ ( 1 - \sum_{a \in A}~ \lambda_a \psi_a ) U_{-} +
                         ( \sum_{a \in A}~ \lambda_a \psi_a ) U_{+} $$

\medskip
These two functions $U_{-}$ and $U_{+}$ satisfy on $X \cap B ( x, \delta )$ the inequalities

$$ T ( . , D ) U_{-} ~\leq~ f,~~~~~ f ~\leq~ T ( . , D ) U_{+} $$ \\ \\

{\bf Recent Extensions}. An analysis of the above results suggests the following more general
one which is obtainable in a rather direct way. \\ \\

{\bf Proposition 2} \\

Suppose $X \subseteq {\bf R}^n$ is a nonvoid, possibly unbounded open subset, and the
continuous nonlinear PDEs (2.1) satisfies the range condition, see (2.7)

\bigskip
(5.12) \quad $ f(x) \in R_x,~~ x \in X_{sol} $

\medskip
where $X_{sol} \subseteq X$ is a nonvoid subset of $X$. Then

\bigskip
(5.13) \quad $ \begin{array}{l}
                  \forall~ A \subseteq X_{sol},~ A ~\mbox{discrete in}~ X ~: \\
                  \\
                  \exists~ U \in {\cal C}^\infty (X) ~: \\
                  \\
                  ~~~ T(X,D) U(x) = f(x),~~ x \in A
               \end{array} $

\bigskip
{\bf Proof} \\
Since $A$ is discrete in $X$, there exists a family $\psi_a \in {\cal C}^\infty (X)$, with
$a \in A$, such that for each given $a \in A$ we have $\psi_a = 1$ on a neighbourhood of $a$,
while for each $b \in A,~ b \neq a$, we have $~\mbox{supp}~ \psi_a ~\cap~ \mbox{supp}~ \psi_b
= \phi$. \\

On the other hand $A \subseteq X_{sol}$ and (5.12) give for each $a \in A$ some real numbers
$\xi_{a, p} \in {\bf R}$, with $p \in {\bf N}^n,~ |p| \leq m$, such that

$$ F(a, \xi_{a, 0}, ~.~.~.~ , \xi_{a, p}, ~.~.~.~ ) = f(a) $$

\medskip
But for each $a \in A$, we can find a polynomial $P_a$ in $x \in {\bf R}^n$ such that

$$ D^p P_a (a) = \xi_{a, p},~~ p \in {\bf N}^n,~ |p| \leq m $$

\medskip
Now we define the function $U \in {\cal C}^\infty (X)$ by

$$ U (x) ~=~ \sum_{a \in A}~ \psi_a (x) P_a (x),~~ x \in X $$

\medskip
and obtain the result in (5.13). \\ \\

{\bf Remark 3} \\

We note that the result in Proposition 2 does not actually need the continuity of the PDE in
(2.1), that is, of $F$ and $f$.

\hfill $\diamondsuit \diamondsuit \diamondsuit$ \\

And now we turn to the basic result in Corollary 2. Here the setup will revert to the case of
smooth nonlinear PDEs in (2.1), (2.2). \\
On the other hand, we shall have to assume certain appropriate stronger alternatives of  the
rather weak range condition (2.7) , alternatives which are defined now. \\
First, it will be convenient to write the smooth nonlinear PDEs in (2.1) in the following
equivalent form

\bigskip
(5.14) \quad $ T (x,D) U (x) ~=~ 0,~~~ x \in X $

\medskip
simply by assimilating the right hand term $f$ in (2.1) into $F$ in (2.2). Now,  for every
$l \in {\bf N}$, let us consider associated with the smooth nonlinear PDE in (5.14) the
following smooth nonlinear one, given by

\bigskip
(5.15) \quad $ T^l (x,D) U (x) = 0,~~~ x \in X $

\medskip
where the partial differential operator $T^l (x,D)$ is defined by

\bigskip
(5.16) \quad $ T^l (x,D) U (x) =
         \sum_{p \in {\bf N}^n,~ |p| \leq l}~ ( D^p ( T (x,D) U (x)  ) )^2,~~~ x \in X $

\medskip
Clearly the smooth nonlinear PDE in (5.15) is equivalent with the system of smooth PDEs

\bigskip
(5.17) \quad $ D^p (T (x,D) U (x) ) = 0,~~~ p \in {\bf N}^n,~ |p| \leq l,~ x \in X $

\medskip
which contains (5.14) among them, for $p = (0, ~.~.~.~ , 0) \in {\bf N}^n$. Therefore the
smooth nonlinear PDE in (5.15) can be of order up to, and including $m + l$. \\
We note that in the case of the smooth PDE in (5.14), the range condition (2.7) takes the
simple and equivalent form

\bigskip
(5.18) \quad $ 0 \in R_x,~~~ x \in X_{sol} $

\medskip
This simple form will be convenient next, when we have to deal with the range condition in the
context of the rather involved smooth PDEs in (5.15), that is, of the systems of smooth
nonlinear PDEs in (5.17). \\
Let us now, according to (2.3), associate with each smooth PDE in (5.15) and at each
$x \in X$, its corresponding range, namely

\bigskip
(5.19) \quad $ R^l_x \subseteq {\bf R},~~~ l \in {\bf N} $

\medskip
For the smooth nonlinear PDE in (5.14) we shall consider now the following stronger version of
the range condition (2.7), namely

\bigskip
(5.20) \quad $ 0 \in R^l_x,~~~ l \in {\bf N},~ x \in X_{sol} $

\medskip
Remarks on the range condition (5.20) which help in clarifying its nature will be given in
section 7. \\ \\

{\bf Corollary 2} \\

Suppose $X \subseteq {\bf R}^n$ is a nonvoid, possibly unbounded open subset, and the smooth
nonlinear PDEs (5.14) satisfies the range condition (5.20) on a nonvoid and at least countable
subset $X_{sol}$ of $X$. \\

Further, suppose given any sequence of integers $l_\nu \in {\bf N}$, with $\nu \in {\bf N}$,
such that $\lim_{~\nu \to \infty}~ l_\nu = \infty$. \\

Then for every countable subset $Z = \{~ z_\nu ~|~ \nu \in {\bf N} ~\} \subseteq X_{sol}$ one
can construct a sequence of smooth functions $s = ( s_\nu ~|~ \nu \in {\bf N} ) \in
( {\cal C}^\infty (X) )^{\bf N}$, such that the {\it error} sequence $ w =
( w_\nu ~|~ \nu \in {\bf N} ) \in ( {\cal C}^\infty (X) )^{\bf N}$ which corresponds to $s$
according to the PDE in (5.14), namely

\bigskip
(5.21) \quad $ w = T (x,D) s $

\medskip
that is, with $w_\nu = T(x,D) s_\nu $, for $\nu \in {\bf N}$, has the {\it vanishing} property

\bigskip
(5.22) \quad $ \begin{array}{l}
                     \forall~ \nu \in {\bf N},~ p \in {\bf N}^n,~ |p| \leq l_\nu ~: \\
                     \\
                     ~~~ D^p w_\nu (z_\nu) = 0
                  \end{array} $

\bigskip
{\bf Proof} \\
Let us take $Z = \{ z_\nu ~|~ \nu \in {\bf N} \}$ a countable subset of $X_{sol}$. Then for any $\nu \in {\bf N}$ we apply
Proposition 1.2 to the discrete subset $A = \{~ z_0, ~.~.~.~ , z_\nu ~\}$ and the smooth PDE in (5.15) which corresponds
to $l_\nu$. \\ \\

{\bf 6. Constructing Solution Algebras with Dense Singularities} \\

The algebras needed in Theorem 1 are constructed as follows. Let again be given any possibly
unbounded nonvoid open subset $X \subseteq {\bf R}^n$. \\
Given any dense subset $X_{sol}$ in $X$, then clearly $X_{sol}$ is at least countable, as
required in Corollary 2. Further, let $l_\nu \in {\bf N}$, with $\nu \in {\bf N}$, be such
that $\lim_{\nu \to \infty}~ l_\nu = \infty$. \\
Let us now take any countable subset $Z = \{~ z_\nu ~|~ \nu \in {\bf N} ~\} \subseteq
X_{sol}$ which is still dense in $X$. \\
Then Corollary 2 yields a sequence of smooth functions

\bigskip
(6.1) \quad $ s = ( s_\nu ~|~ \nu \in {\bf N} ) \in ( {\cal C}^\infty (X) )^{\bf N} $

\medskip
with the following {\it error} property regarding the smooth nonlinear PDE in (2.1). Let us
denote

\bigskip
(6.2) \quad $ w = T(x,D) s - u(f) $

\medskip
that is, $w = ( w_\nu ~|~ \nu \in {\bf N} ) \in ( {\cal C}^\infty (X) )^{\bf N}$, where
$w_\nu = T(x,D) s_\nu - f$, with $\nu \in {\bf N}$. Then in view of (1.3.8), the {\it
vanishing} property (5.22) means that

\bigskip
(6.3) \quad $ w \in {\cal J}_\Sigma (X) $

\medskip
where $\Sigma = X \setminus Z$ will clearly satisfy (3.3), and in particular, the requirements
in Theorem 1. Here we note that, in more detail, the vanishing property (5.22) means

\bigskip
(6.4) \quad $ \begin{array}{l}
                \forall~ \nu \in {\bf N},~ p \in {\bf N}^n,~ |p| \leq l_\nu ~: \\
                              \\
                              ~~~D^p ~ (~ T(x,D) s_\nu - f ~) (z_\nu ) = 0
               \end{array} $

\medskip
Now if we take the generalized function $U$ in (4.1), then in the sense of (3.19), we shall
clearly have satisfied (4.2). \\

In the above construction of the generalized solutions $U$ through the respective sequences of
smooth functions $s$, we have had the liberty, under the conditions of Theorem 1, to choose a
large variety of sequences of integers $l_\nu$, as well as countable dense subsets $Z$. This
means that, correspondingly, the above construction can deliver a large variety of generalized
solutions $U$. \\ \\

{\bf 7. Remarks of the Range Condition} \\

Here we present a certain clarification of the meaning of the range condition (5.20) used in
the vanishing result in Corollary 2, and thus, in the basic existence of solutions result in
Theorem 1. For that purpose one can start with the equivalence, for each given $l \in
{\bf N}$, between the smooth PDE in (5.15), and on the other hand, the system of smooth PDEs
in (5.17). In view of this equivalence, it is natural to relate the ranges $R^l_x$ in (5.20),
which correspond to the PDE in (5.15), with the ranges which correspond to each individual
PDE in the system (5.17). \\
Let us therefore, for $p \in {\bf N}^n$, denote by $R_{p,x}$, with $x \in X$, the ranges of
the corresponding PDE in the system (5.17), see (2.3). \\
In this respect, for any $l \in {\bf N}$ and nonvoid subset $X_{sol} \subseteq X$, we have in
view of (5.16) the implication

\bigskip
(7.1) \quad $ \begin{array}{l}
               \left(~~ 0 \in R^l_x,~~~ x \in X_{sol} ~~\right) ~~~\Longrightarrow \\
                                   \\
               ~~~~~~~~~~~~~~\Longrightarrow~~~ \left(~~ 0 \in R_{p,x},~~~ p \in
                        {\bf N}^n,~ |p| \leq l,~ x \in X_{sol} ~~\right)
                          \end{array} $

\medskip
In this way, we obtain the {\it necessary} condition for the range condition (5.20), namely

\bigskip
(7.2) \quad $ ~~~~~~~~~~(5.20) ~~~\Longrightarrow~~~ \left(~~ 0 \in
                    R_{p,x},~~~ p \in {\bf N}^n,~ x \in X_{sol} ~~\right) $

\medskip
and clearly, the ranges $R_{p,x}$ are easier to compute than the ranges $R^l_x$ which appear
in (5.20), since for any given $l \in {\bf N}$, each of the PDEs in the system (5.17) is
significantly simpler than the PDE in (5.15). \\

In view of the equivalence between (5.15) and (5.17), the question arises to what extent does
the converse implication hold in (7.2), which would then obviously give a natural {\it
characterization} of the range condition (5.20). \\
Clearly, this converse implication in (7.2) would result, if we can establish the easier
converse implication in (7.1), for each particular $l \in {\bf N}$. We proceed now to clarify
this issues. \\

Before we go further, it is useful to note the following property of the ranges $R_x$, with
$x \in X$, in the general definition (2.3). Given the smooth nonlinear PDE in (2.1), (2.2), we
shall consider it in the form (5.14). We recall that earlier we denoted by $m^*$ the number of
real variable arguments in $F$ minus $n$. Then we clearly have

\bigskip
(7.3) \quad $ R_x ~=~ F ( \{ x \} \times {\bf R}^{m^*} ),~~~ x \in X $

\medskip
Now returning to (5.17), obviously this smooth nonlinear PDE is of the form

\bigskip
(7.4) \quad $ F_p ( x, U(x),~.~.~.~, D^r U(x),~.~.~.~ ) ~=~ 0,~~~ x \in X $

\medskip
where $r \in {\bf N}^n,~ |r| \leq m + |p|$. Thus

\bigskip
(7.5) \quad $ R_{p,x} = F_p ( \{ x \} \times {\bf R}^{ ( m + |p| )^*} ),~~~ x \in X $

\medskip
On the other hand, the smooth nonlinear PDE in (5.15) is of the form

\bigskip
(7.6) \quad $ F^l (x, U(x),~.~.~.~, D^r U(x),~.~.~.~ ) ~=~ 0,~~~ x \in X $

\medskip
with $r \in {\bf N}^n,~ |r| \leq m + l$. Hence

\bigskip
(7.7) \quad $ R^l_x = F^l ( \{ x \} \times {\bf R}^{(m + l)^*} ),~~~ x \in X $

\medskip
In connection with the sought after converse of the implication (7.1) let us fix $l \in
{\bf N}$. It will be convenient to consider the functions $F_p$ in (7.4), with $p \in
{\bf N}^n,~ |p| \leq l$, as being defined on the same domain $X \times {\bf R}^{(m+l)^*}$ as
the function $F^l$ in (7.6). This extension is done in the obvious natural way illustrated in
the particular case in (7.15). Then the converse implication in (7.1) can be formulated as

\bigskip
(7.8) \quad $ \begin{array}{l}
                             \left(~ 0 \in \bigcap_{p \in {\bf N}^n,~ |p| \leq l}~ F_p ( \{ x \} \times
                                                    {\bf R}^{(m+l)^*} ),~~~ x \in X_{sol} ~\right) ~~~\Longrightarrow~~~ \\
                             \\
                              ~~~~~~~~~~~~~~~~~~\Longrightarrow~~~ \left(~ 0 \in F^l ( \{ x \} \times {\bf R}^{(m+l)^*} ),~~~ x \in X_{sol} ~\right)
                          \end{array} $

\medskip
This is of course the same with the requirement that the condition

\bigskip
(7.9) \quad $ \begin{array}{l}
                    \forall~ x \in X_{sol},~ p \in {\bf N}^n,~ |p| \leq l ~: \\
                               \\
                    \exists~ \xi_p = ( \xi_{p,1},~.~.~.~, \xi_{p,(m+l)^*} ) \in
                           {\bf R}^{(m+l)^*} ~: \\
                               \\
                                ~~~ F_p (x, \xi_p) = 0
               \end{array} $

\medskip
imply the relation

\bigskip
(7.10) \quad $ \begin{array}{l}
                  \forall~ x \in X_{sol} ~: \\
                             \\
                  \exists~ \eta^l = ( \eta_1,~.~.~.~, \eta_{(m+l)^*} ) \in {\bf R}^{(m+l)^*} ~: \\
                             \\
                      ~~~F^l (x, \eta^l) = 0
               \end{array} $

\medskip
We also note that (5.16) is the same with the relation

\bigskip
(7.11) \quad $ F^l ~=~ \sum_{p \in {\bf N}^n,~ |p| \leq l}~ ( F_p )^2 ~~~
                      \mbox{on}~~ X \times {\bf R}^{(m+l)^*} $

\medskip
In other words, the converse implication in (7.1) means that, given $x \in X_{sol}$, the
existence of zeros for each of the right hand terms in (7.11) must imply the existence of a
{\it common} zero for all these terms, that is, the existence of a zero for the left hand
term. \\
It will be useful to illustrate the above in some simple cases, in order to gain some
understanding about the possibility of this converse implication. \\

{\bf The Smooth Linear Case}. For convenience we start with the case of smooth linear PDEs
(2.1) of the general form

\bigskip
(7.12) \quad $ \sum_{q \in {\bf N}^n,~|q| \leq m}~ c_q (x) D^q U (x) ~=~ f (x),~~~ x \in X $

\medskip
with $c_q, f \in {\cal C}^\infty (X)$. We shall consider (7.12) in terms of (5.14), and our
interest is in the respective nonlinear smooth PDE in (5.15), where $l \in {\bf N}$ is fixed
for the moment, namely

\bigskip
(7.13) \quad $ \begin{array}{l}
                \sum_{p \in {\bf N}^n,~|p| \leq l}~
                 ( \sum_{q \in {\bf N}^n,~|q| \leq m}~ D^p ( c_q (x) D^q U (x) - f (x) ) )^2
                              = 0, \\
                              ~~~~~~~~~~~~~~~~~~~~~~~~~~~~~~~~~~~~~~~~~~~~~~~~~~~~~~~~~~~~~~~~~~~~~~~ x \in X
                 \end{array} $

\medskip
Then for a given $p \in {\bf N}^n$, the corresponding PDE in the system (5.17) becomes

\bigskip
(7.14) \quad $ \sum_{q \in {\bf N}^n,~|q| \leq m}~ D^p ( c_q (x) D^q U (x) - f (x) ) ~=~ 0,
                                            ~~~ x \in X $

\medskip
Given $l \in {\bf N}$ fixed, then in view of (7.4), (7.14), it is obvious that $F_p$, with
$p \in {\bf N}^n,~ |p| \leq l$, will be of the form

\bigskip
(7.15) \quad $ \begin{array}{l}
                   F_p (x, \xi_1,~.~.~.~, \xi_{(m+l)^*}) = d_p(x) +
                           \sum_{1 \leq i \leq {(m+|p|)^*}}~ d_{p,i}(x) \xi_i, \\
                               ~~~~~~~~~~~~~~~~~~~~~~~~~~~~~~~~~~~~~~~~~~~~~~~~~~~~~~~~
                                            x \in X,~ \xi_i \in {\bf R}
               \end{array} $

\medskip
with $d_p,~d_{p,i} \in {\cal C}^\infty (X)$. Consequently (7.11) gives

\bigskip
(7.16) \quad $ \begin{array}{l}
                    F^l (x, \xi_1,~.~.~.~, \xi_{(m+l)^*}) = \\
                     = \sum_{p \in {\bf N}^n,~ |p| \leq l}~(~d_p(x) +
                               \sum_{1 \leq i \leq {(m+|p|)^*}}~
                                              d_{p,i}(x) \xi_i ~)^2, \\
                             ~~~~~~~~~~~~~~~~~~~~~~~~~~~~~~~~~~~~~~~~~~~~~~~~~~~~~~~
                                    x \in X,~ \xi_i \in {\bf R}
                \end{array} $

\medskip
Clearly, for every given $x \in X$, each $F_p$ in (7.15) vanishes on a hyperplane of dimension
at least $(m+|p|)^*-1$ in ${\bf R}^{(m+l)^*}$. Therefore, the vanishing of $F^l$, which in
this case is equivalent to the converse implication in (7.1), is the same with the existence
of a solution in $\xi_i$ of the linear system

\bigskip
(7.17) \quad $ \sum_{1 \leq i \leq {(m+|p|)^*}}~ d_{p,i}(x) \xi_i =
                                      - d_p(x),~~~ p \in {\bf N}^n,~ |p| \leq l $

\medskip
of $(m+l)^*$ equations with the {\it same} number of unknowns $\xi_i$. Let us for each
$x \in X$ denote by $P^l(x)$ and $Q^l(x)$, respectively, the coefficient matrix of the
unknowns $\xi_i$ and the matrix of all coefficients in (7.17). Then the {\it converse}
implication, hence the {\it equivalence} in (7.1), takes the form

\bigskip
(7.18) \quad $ \mbox{rank}~ P^l(x) ~=~ \mbox{rank}~ Q^l(x),~~~ x \in X_{sol} $

\medskip
Finally, in view of (7.2), we obtain the {\it equivalent} form of the range condition (5.20)
as given by

\bigskip
(7.19) \quad $ \mbox{rank}~ P^l(x) ~=~ \mbox{rank}~ Q^l(x),~~~ x \in X_{sol},~ l \in {\bf N} $

\medskip
Without getting involved in more detailed technical arguments, it is easy to see that for each
given $l \in {\bf N}$ and $x \in X$, the linear system (7.17) is of lower diagonal block type.
And the diagonal blocks are made up, for each $p \in {\bf N}^n$, of those coefficients
$d_{p,i}(x)$ which correspond to $1 \leq i \leq m + |p|$ for which the respective terms belong
to the principal part of the smooth linear PDE in (7.14). It follows that for nondegenerate
smooth linear PDEs in (7.12), it is likely to have on a dense subset $X_{sol}$ of $X$
satisfied the condition (7.18) in the stronger form

\bigskip
(7.20) \quad $ \mbox{rank}~ P^l(x) ~=~ \mbox{rank}~ Q^l(x) ~=~ (m + l)^*,~~~ x \in X_{sol} $

\medskip
for each $l \in {\bf N}$. \\
Therefore, Theorem 1 in section 4 does to that extent overcome the Lewy impossibility for all
such smooth linear PDEs in (7.12). \\

{\bf Open Problem}. One can try to elaborate the full details of the above diagonal type
argument which may lead to the validity of (7.18) or (7.20) on a dense subset $X_{sol}$ of $X$
and for all $l \in {\bf N}$. \\ \\

{\bf 8. Comments} \\

1. In 1957, H. Lewy showed that the surprisingly simple smooth linear PDE

\bigskip
(8.1) \quad $ \begin{array}{l}
                     ( D_x + i D_y - 2 (x + i y) D_z ) U(x, y, z) ~=~ f(x, y, z), \\
                                  ~~~~~~~~~~~~~~~~~~~~~~~~~~~~~~~~~~~~~~~~~~~~~~~~~~~
                                   (x, y, z) \in {\bf R}^3
              \end{array} $

\medskip
does not have any distribution solution $U \in {\cal D}^\prime$ in any neighbourhood of any
point of ${\bf R}^3$, for a large class of right hand terms $f \in {\cal C}^\infty({\bf R}^3)$.
Further studies of this {\it impossibility} to solve smooth linear PDEs in Schwartz
distributions followed, and a similar impossibility has been proved to hold also with respect
to the Sato hyperfunctions, see details and references in Rosinger [7, pp. 37-39]. \\

2. The first time the Lewy impossibility was overcome was in 1994, in Oberguggenberger \&
Rosinger, where arbitrary continuous nonlinear PDEs were given solutions which can be
assimilated with usual measurable functions. In this way, the mentioned result overcomes the
Lewy impossibility with a very large margin. \\
Recently, that result has been improved in Anguelov \& Rosinger by showing that the respective
solutions can be assimilated with Hausdorff continuous functions. \\

3. The main aim of the method in Oberguggenberger \& Rosinger and also in Anguelov \& Rosinger
was to apply some of the most basic possible mathematical constructions - such as order
completion - to the solution of the largest possible classes of nonlinear PDEs. \\
One of the first result in this regard is that in Proposition 1 in section 5, conjectured and
proved in August 1990 by M. Oberguggenberger in the case of bounded domains $X$. This result
underlies some of the first general existence of solutions results, and as one of its many
byproducts, it leads to the first ever overcoming of the Lewy impossibility. See for further
details Oberguggenberger \& Rosinger [pp. VII-XII, 3-10]. \\

4. The systematic development of the algebraic-differential nonlinear method for global
generalized solutions of large classes of nonlinear PDEs started back in the 1960s in Rosinger
[1,2] and its development continued in Rosinger [3-17], Rosinger \& Walus [1,2], Mallios \&
Rosinger [1,2], Colombeau [1,2], Biagioni, Oberguggenberger. \\
Since the mid 1980s, several dozen other publications by a variety of authors appeared. \\

During October-December 1997, the Erwin Schroedinger International Institute for Mathematical
Physics of the University of Vienna, Austria, organized the first Workshop on the Nonlinear
Theory of Generalized Functions, the proceedings of which, published in 1999,  see Grosser \&
Hoermann \& Kunzinger \& Oberguggenberger, contain more than two dozen of the recent
contributions, together with relevant references. \\
A more recent development is presented in Grosser et.al. \\

In April 2000, the ICGF2000 - International Conference on Generalized Functions, Linear and
Nonlinear Problems was held at the Universite des Antilles at de la Guyane, Guadeloupe, where
more than fifty papers on a large variety of subjects in the field were presented. \\

6. In 1999, the field at large to which its subject - together with that of this paper -
belongs, namely, the {\it nonlinear algebraic differential theory of generalized functions},
with the respective {\it differential algebras of generalized functions}, was included by the
American Mathematical Society in their AMS Classification 2000, under the heading : \\

46F30 Generalized functions for nonlinear analysis \\

See further details at : ~~ www.ams.org/index/msc/46Fxx.html \\

Earlier, the subject of the same heading 46F30 used to be 'Distributions, generalized
functions, distribution spaces', which was now changed, as perhaps it was felt that it was
sounding not enough nonlinear ... \\

\end{document}